\documentclass[12pt]{article}

\usepackage[geometry]{ifsym}
\usepackage{verbatim}
\usepackage{fancyhdr}
\usepackage{graphicx}
\usepackage{mathrsfs}
\usepackage{amssymb}
\usepackage{cite}
\usepackage{epsfig}
\usepackage{pst-poly}     
\usepackage{pst-all}      

\usepackage{amsfonts,amssymb, amsmath, latexsym}

\newtheorem{cor}{Corollary}[section]
\newtheorem{lem}{Lemma}[section]
\newtheorem{thm}{Theorem}[section]

\newenvironment{pf}[1][Proof]{\noindent\textbf{#1.} }{\hfill\rule{1mm}{2mm}}

\makeatletter \@addtoreset{equation}{section} \makeatother

\begin{document}

\title{Total and paired domination numbers of toroidal meshes~\thanks {The work was supported by NNSF
of China (No. 11071233).}}
\author
{Fu-Tao Hu,\quad Jun-Ming Xu\footnote{Corresponding
author:
xujm@ustc.edu.cn}\ \\ \\
{\small Department of Mathematics}  \\
{\small University of Science and Technology of China}\\
{\small Hefei, Anhui, 230026, China} }
\date{}
\maketitle

\begin{quotation}
\textbf{Abstract}: Let $G$ be a graph without isolated vertices. The
total domination number of $G$ is the minimum number of vertices
that can dominate all vertices in $G$, and the paired domination
number of $G$ is the minimum number of vertices in a dominating set
whose induced subgraph contains a perfect matching. This paper
determines the total domination number and the paired domination
number of the toroidal meshes, i.e., the Cartesian product of two
cycles $C_n$ and $C_m$ for any $n\ge 3$ and $m\in\{3,4\}$, and gives
some upper bounds for $n, m\ge 5$.



\end{quotation}

\section{Introduction}

For notation and graph-theoretical terminology not defined here we
follow \cite{x03}. Specifically, let $G=(V,E)$ be an undirected
graph without loops, multi-edges and isolated vertices, where
$V=V(G)$ is the vertex-set and $E=E(G)$ is the edge-set, which is a
subset of $\{xy|\ xy$ is an unordered pair of $V \}$. A graph $G$ is
{\it nonempty} if $E(G)\ne \emptyset$. Two vertices $x$ and $y$ are
{\it adjacent} if $xy\in E(G)$. For a vertex $x$, denote $N(x)=\{y:
xy\in E(G)\}$ be the {\it neighborhood} of $x$. For a subset
$D\subseteq V(G)$, we use $G[D]$ to denote the subgraph of $G$
induced by $D$. We use $C_n$ and $P_n$ to denote a cycle and a path
of order $n$, respectively, throughout this paper.

A subset $D\subseteq V(G)$ is called a {\it dominating set} if
$N(x)\cap D\ne \emptyset$ for each vertex $x\in V(G)\setminus D$.
The {\it domination number} $\gamma(G)$ is the minimum cardinality
of a dominating set. A thorough study of domination appears in
\cite{hhs98a, hhs98b}. A subset $D\subseteq V(G)$ of $G$ is called a
{\it total dominating set}, introduced by Cockayne {\it et
al.}~\cite{cdh80}, if $N(x)\cap D\ne \emptyset$ for each vertex
$x\in V(G)$ and the {\it total domination number} of $G$, denoted by
$\gamma_t(G)$, is the minimum cardinality of a total dominating set
of $G$. The total domination in graphs has been extensively studied
in the literature. A survey of selected recent results on this topic
is given in~\cite{h09} by Henning.

A dominating set $D$ of $G$ is called to be {\it paired}, introduced
by Haynes and Slater ~\cite{hs95,hs98}, if the induced subgraph
$G[D]$ contains a perfect matching. The {\it paired domination
number} of $G$, denoted by $\gamma_p(G)$, is the minimum cardinality
of a paired dominating set of $G$. Clearly,
 $\gamma(G)\le\gamma_t(G)\le \gamma_p(G)$
since a paired dominating set is also a total dominating set of $G$,
and $\gamma_p(G)$ is even. Pfaff, Laskar and Hedetniemi~\cite{plh83}
and Haynes and Slater~\cite{hs98} showed that the problems
determining the total-domination and the paired-domination for
general graphs are NP-complete. Some exact values of
total-domination numbers and paired-domination numbers for some
special classes of graphs have been determined by several authors.
In particularly, $\gamma_t(P_n \times P_m)$ and $\gamma_p(P_n \times
P_m)$ for $2\le m\le 4$ are determined by Gravier~\cite{g02}, and
Proffitt, Haynes and Slater~\cite{phs01}, respectively.

Use $G_{n,m}$ to denote the toroidal meshes, i.e., the Cartesian
product $C_n \times C_m$ of two cycles $C_n$ and $C_m$. Klav\v{z}ar
and Seifter~\cite{ks95} determined $\gamma(G_{n,m})$ for any $n\ge
3$ and $m\in\{3,4,5\}$. In this paper, we obtain the following
results.

$$
\begin{array}{rl}
&\gamma_{t}(G_{n,3})=\lceil \frac{4n}{5}\rceil;\\
&\gamma_{p}(G_{n,3})=\left\{
 \begin{array}{ll}
\lceil\frac{4n}{5}\rceil  & {\rm if}\ n\equiv 0,2,4\,({\rm mod}\,5),\\
\lceil\frac{4n}{5}\rceil+1  & {\rm if}\ n\equiv 1,3\,({\rm mod}\,5);\\
 \end{array}\right.\\
&\gamma_{t}(G_{n,4})=\gamma_{p}(G_{n,4})=\left\{
 \begin{array}{ll}
n &  {\rm if}\ n\equiv 0\,({\rm mod}\,4),\\
n+1 & {\rm if}\ n\equiv 1,3\,({\rm mod}\,4),\\
n+2 & {\rm if}\ n\equiv 2\,({\rm mod}\,4).
 \end{array}\right.
\end{array}
 $$


\section{Preliminary results}

In this section, we recall some definitions, notations and results
used in the proofs of our main results. Throughout this paper, we
assume that a cycle $C_n$ has the vertex-set
$V(C_n)=\{1,\ldots,n\}$.

Use $G_{n,m}$ to denote the toroidal meshes, i.e., the Cartesian
product $C_n \times C_m$, which is a graph with vertex-set
$V(G_{n,m})=\{x_{ij}|\ 1\leq i \leq n, 1\leq j \leq m\}$ and two
vertices $x_{ij}$ and $x_{i'j\,'}$ being linked by an edge if and
only if either $i=i'\in V(C_n)$ and $jj\,'\in E(C_m)$, or $j=j\,'\in
V(C_m)$ and $ii'\in E(C_n)$.

Let $Y_i=\{x_{ij}|\ 1\leq j \leq m\}$ for $1\leq i\leq n$, called a
set of {\it vertical vertices} in $G_{n,m}$.

In~\cite{gs02}, Gavlas and Schultz defined an efficient total
dominating set, which is such a total dominating set $D$ of $G$ that
$|N(v)\cap D|=1$ for every $v\in V(G)$. The related research results
can be found in~\cite{ds03, gs02, hx08}.

\begin{lem}\label{lem2.1}{\rm(Gavlas~and~Schult\cite{gs02})}\ If a
graph $G$ has an efficient total dominating set $D$, then the
edge-set of the subgraph $G[D]$ forms a perfect matching, and so the
cardinality of $D$ is even, and $\{N(v): v\in D\}$ partitions
$V(G)$.
\end{lem}

\begin{lem}\label{lem2.2}
Let $G$ be a $k$-regular graph of order $n$. Then $\gamma_t(G)\ge
\frac{n}{k}$, with equality if and only if $G$ has an efficient
total dominating set.
\end{lem}

\begin{pf}
Since $G$ is $k$-regular, each $v\in V(G)$ can dominate at most $k$
vertices. Thus $\gamma_t(G)\ge \frac{n}{k}$. It is easy to observe
that the equality holds if and only if there exists a total
dominating set $D$ such that $\{N(v): v\in D\}$ partitions $V(G)$,
equivalently, $D$ is an efficient total dominating set.
\end{pf}

\begin{lem}\label{lem2.3}
$\gamma_{t}(G_{n,m})=\gamma_{p}(G_{n,m})=\frac{nm}{4}$ for
$n,m\equiv \,0~({\rm mod}\,4)$.
\end{lem}

\begin{pf}
Let $D=\{x_{ij},x_{i(j+1)},x_{(i+2)(j+2)},x_{(i+2)(j+3)}: ~i,j\equiv
\,1~({\rm mod}\,4)\}$, where $1\le i\le n$ and $1\le j\le m$.
Figure~\ref{f1} is such a set $D$ in $G_{8,4}$. It is easy to see
that $D$ is a paired dominating set of $G_{n,m}$ with cardinality
$\frac{nm}{4}$. Thus, $\gamma_{p}(G_{n,m})\le \frac{nm}{4}$.

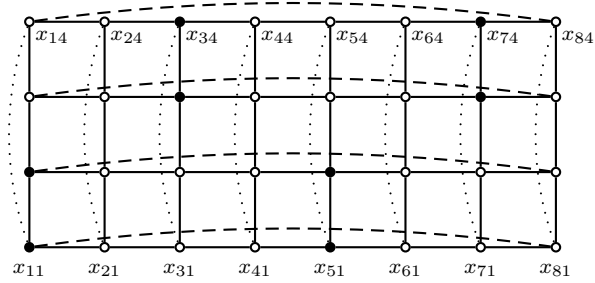
\begin{figure}[h]
\begin{center}
\begin{pspicture}(0,.5)(8.5,5)

\cnode*(1,1){2pt}{11}\rput(1,.7){\scriptsize $x_{11}$}
\cnode*(1,2){2pt}{12}  \cnode(1,3){2pt}{13}
\cnode(1,4){2pt}{14}\rput(1.3,3.8){\scriptsize $x_{14}$}

\cnode(2,1){2pt}{21}\rput(2,.7){\scriptsize $x_{21}$}
\cnode(2,2){2pt}{22}  \cnode(2,3){2pt}{23}
\cnode(2,4){2pt}{24}\rput(2.3,3.8){\scriptsize $x_{24}$}

\cnode(3,1){2pt}{31}\rput(3,.7){\scriptsize $x_{31}$}
\cnode(3,2){2pt}{32}  \cnode*(3,3){2pt}{33}
\cnode*(3,4){2pt}{34}\rput(3.3,3.8){\scriptsize $x_{34}$}

\cnode(4,1){2pt}{41}\rput(4,.7){\scriptsize $x_{41}$}
\cnode(4,2){2pt}{42}  \cnode(4,3){2pt}{43}
\cnode(4,4){2pt}{44}\rput(4.3,3.8){\scriptsize $x_{44}$}

\cnode*(5,1){2pt}{51}\rput(5,.7){\scriptsize $x_{51}$}
\cnode*(5,2){2pt}{52}  \cnode(5,3){2pt}{53}
\cnode(5,4){2pt}{54}\rput(5.3,3.8){\scriptsize $x_{54}$}

\cnode(6,1){2pt}{61}\rput(6,.7){\scriptsize $x_{61}$}
\cnode(6,2){2pt}{62}  \cnode(6,3){2pt}{63}
\cnode(6,4){2pt}{64}\rput(6.3,3.8){\scriptsize $x_{64}$}
\cnode(7,1){2pt}{71}\rput(7,.7){\scriptsize $x_{71}$}
\cnode(7,2){2pt}{72}  \cnode*(7,3){2pt}{73}
\cnode*(7,4){2pt}{74}\rput(7.3,3.8){\scriptsize $x_{74}$}

\cnode(8,1){2pt}{81}\rput(8,.7){\scriptsize $x_{81}$}
\cnode(8,2){2pt}{82}  \cnode(8,3){2pt}{83}
\cnode(8,4){2pt}{84}\rput(8.3,3.8){\scriptsize $x_{84}$}

\ncline{11}{12} \ncline{12}{13}  \ncline{11}{21} \ncline{12}{22}
\ncline{13}{23}  \ncline{13}{14}  \ncline{14}{24}
\ncarc[linestyle=dotted,arcangle=20]{11}{14} \ncline{21}{22}
\ncline{22}{23}  \ncline{21}{31} \ncline{22}{32}  \ncline{23}{33}
\ncline{23}{24}  \ncline{24}{34}
\ncarc[linestyle=dotted,arcangle=20]{21}{24} \ncline{31}{32}
\ncline{32}{33}  \ncline{31}{41} \ncline{32}{42}  \ncline{33}{43}
\ncline{33}{34}  \ncline{34}{44}
\ncarc[linestyle=dotted,arcangle=20]{31}{34} \ncline{41}{42}
\ncline{42}{43}  \ncline{41}{51} \ncline{42}{52}  \ncline{43}{53}
\ncline{43}{44}  \ncline{44}{54}
\ncarc[linestyle=dotted,arcangle=20]{41}{44} \ncline{51}{52}
\ncline{52}{53}  \ncline{51}{61} \ncline{52}{62}  \ncline{53}{63}
\ncline{53}{54}  \ncline{54}{64}
\ncarc[linestyle=dotted,arcangle=20]{51}{54} \ncline{61}{62}
\ncline{62}{63}  \ncline{61}{71} \ncline{62}{72}  \ncline{63}{73}
\ncline{63}{64}  \ncline{64}{74}
\ncarc[linestyle=dotted,arcangle=20]{61}{64} \ncline{71}{72}
\ncline{72}{73}  \ncline{71}{81} \ncline{72}{82}  \ncline{73}{83}
\ncline{73}{74}  \ncline{74}{84}
\ncarc[linestyle=dotted,arcangle=20]{71}{74} \ncline{81}{82}
\ncline{82}{83}  \ncline{83}{84}
\ncarc[linestyle=dotted,arcangle=20]{81}{84}

\ncarc[linestyle=dashed]{11}{81}  \ncarc[linestyle=dashed]{12}{82}
\ncarc[linestyle=dashed]{13}{83} \ncarc[linestyle=dashed]{14}{84}

\end{pspicture}
\caption{\label{f1}\footnotesize The minimum total (paired)
dominating set (bold vertices) of $G_{8,4}$}
\end{center}
\end{figure}

By Lemma~\ref{lem2.2}, $\gamma_{t}(G_{n,m})\ge \frac{nm}{4}=n$.
Since $\gamma_t(G_{n,m})\le \gamma_{p}(G_{n,m})$,
$\gamma_{t}(G_{n,m})=\gamma_{p}(G_{n,m})=\frac{nm}{4}$.
\end{pf}


\section{Total and paired domination number of $G_{n,3}$}

In this section, we determine the exact values of the total and the
paired domination numbers of $G_{n,3}$, which can be stated the
following theorem.

\begin{thm}\label{thm3.1} For any $n\geq 3$,
$$\gamma_{t}(G_{n,3})=\left\lceil \frac{4n}{5}\right\rceil$$ and
$$\gamma_{p}(G_{n,3})=\left\{
 \begin{array}{ll}
\lceil\frac{4n}{5}\rceil, & {\rm if}\ n\equiv 0,2,4\,({\rm mod}\,5);\\
\lceil\frac{4n}{5}\rceil+1,& {\rm if}\ n\equiv 1,3\,({\rm mod}\,5).
 \end{array}\right.$$
\end{thm}

\begin{pf}
Let $D$ be a minimum total dominating set of $G_{n,3}$. First, we
may assume that $|Y_i\cap D|\le 2$ for any $1\le i\le n$. Indeed, if
$|Y_i\cap D|=3$ for some $i\notin\{1,n\}$, then the set
$D'=(D\setminus \{x_{i1}, x_{i3}\})\cup \{x_{(i-1)2}, x_{(i+1)2}\}$
is also a total dominating set of $G_{n,3}$ with $|D'|=|D|$.

Let $\alpha_k$ be the number of $i$'s for which $|Y_i\cap D|=k$ for
$1\le i\le n$ and $0\le k\le 2$. Then we have
\begin{equation}\label{e3.1}
\alpha_0+\alpha_1+\alpha_2=n.
\end{equation}

Assume $|Y_i\cap D|=0$ for some $i\notin\{1,n\}$. At least one of
$|Y_{i-1}\cap D|$ and $|Y_{i+1}\cap D|$ is 2 since the three
vertices in $Y_i$ should be dominated by $D$, which means that
\begin{equation}\label{e3.2}
2\alpha_2-\alpha_0\ge 0.
\end{equation}

If $|Y_i\cap D|=2$ for some $i$ with $1\le i\le n$, then the two
vertices in $Y_i\cap D$ can dominate at most 7 vertices. Since any
vertex $x\in D$ can dominate at most 4 vertices, we have
\begin{equation}\label{e3.3}
4\alpha_1+7\alpha_2\ge 3n.
\end{equation}

The sum of (\ref{e3.1}), (\ref{e3.2}) and (\ref{e3.3}) implies
$$5\alpha_1+10\alpha_2\ge 4n,$$
and, hence,
 \begin{equation}\label{e3.4}
 \gamma_{t}(G_{n,3})=|D|=\alpha_1+2\alpha_2\ge \left\lceil \frac{4n}{5}\right\rceil.
 \end{equation}

\begin{figure}[ht]
\begin{center}
\begin{pspicture}(0,.5)(11,4)

\cnode(1,1){2pt}{11}\rput(1,.7){\scriptsize $x_{11}$}
\cnode*(1,2){2pt}{12}\rput(1.3,1.8){\scriptsize $x_{12}$}
\cnode(1,3){2pt}{13}\rput(1.3,2.8){\scriptsize $x_{13}$}

\cnode(2,1){2pt}{21}\rput(2,.7){\scriptsize $x_{21}$}
\cnode*(2,2){2pt}{22}\rput(2.3,1.8){\scriptsize $x_{22}$}
\cnode(2,3){2pt}{23}\rput(2.3,2.8){\scriptsize $x_{23}$}

\cnode(3,1){2pt}{31}\rput(3,.7){\scriptsize $x_{31}$}
\cnode(3,2){2pt}{32}\rput(3.3,1.8){\scriptsize $x_{32}$}
\cnode(3,3){2pt}{33}\rput(3.3,2.8){\scriptsize $x_{33}$}

\cnode*(4,1){2pt}{41}\rput(4,0.7){\scriptsize $x_{41}$}
\cnode(4,2){2pt}{42}\rput(4.3,1.8){\scriptsize $x_{42}$}
\cnode*(4,3){2pt}{43}\rput(4.3,2.8){\scriptsize $x_{43}$}

\cnode(5,1){2pt}{51}\rput(5,0.7){\scriptsize $x_{51}$}
\cnode(5,2){2pt}{52}\rput(5.3,1.8){\scriptsize $x_{52}$}
\cnode(5,3){2pt}{53}\rput(5.3,2.8){\scriptsize $x_{53}$}

\cnode(6,1){2pt}{61}\rput(6,0.7){\scriptsize $x_{61}$}
\cnode*(6,2){2pt}{62}\rput(6.3,1.8){\scriptsize $x_{62}$}
\cnode(6,3){2pt}{63}\rput(6.3,2.8){\scriptsize $x_{63}$}

\cnode(7,1){2pt}{71}\rput(7,0.7){\scriptsize $x_{71}$}
\cnode*(7,2){2pt}{72}\rput(7.3,1.8){\scriptsize $x_{72}$}
\cnode(7,3){2pt}{73}\rput(7.3,2.8){\scriptsize $x_{73}$}

\cnode(8,1){2pt}{81}\rput(8,.7){\scriptsize $x_{81}$}
\cnode(8,2){2pt}{82}\rput(8.3,1.8){\scriptsize $x_{82}$}
\cnode(8,3){2pt}{83}\rput(8.3,2.8){\scriptsize $x_{83}$}

\cnode*(9,1){2pt}{91}\rput(9,.7){\scriptsize $x_{91}$}
\cnode(9,2){2pt}{92}\rput(9.3,1.8){\scriptsize $x_{92}$}
\cnode*(9,3){2pt}{93}\rput(9.3,2.8){\scriptsize $x_{93}$}

\cnode(10,1){2pt}{101}\rput(10,0.7){\scriptsize $x_{(10)1}$}
\cnode(10,2){2pt}{102}\rput(10.5,1.8){\scriptsize $x_{(10)2}$}
\cnode(10,3){2pt}{103}\rput(10.5,2.8){\scriptsize $x_{(10)3}$}

\ncline{11}{12} \ncline{12}{13} \ncarc[linestyle=dotted,arcangle=20]{11}{13} \ncline{11}{21} \ncline{12}{22}  \ncline{13}{23}
\ncline{21}{22} \ncline{22}{23} \ncarc[linestyle=dotted,arcangle=20]{21}{23}  \ncline{21}{31} \ncline{22}{32} \ncline{23}{33}
\ncline{31}{32} \ncline{32}{33} \ncarc[linestyle=dotted,arcangle=20]{31}{33} \ncline{31}{41} \ncline{32}{42}  \ncline{33}{43}
\ncline{41}{42} \ncline{42}{43} \ncarc[linestyle=dotted,arcangle=20]{41}{43} \ncline{41}{51} \ncline{42}{52}  \ncline{43}{53}
\ncline{51}{52} \ncline{52}{53} \ncarc[linestyle=dotted,arcangle=20]{51}{53} \ncline{51}{61} \ncline{52}{62}  \ncline{53}{63}
\ncline{61}{62} \ncline{62}{63} \ncarc[linestyle=dotted,arcangle=20]{61}{63} \ncline{61}{71} \ncline{62}{72}  \ncline{63}{73}
\ncline{71}{72} \ncline{72}{73} \ncarc[linestyle=dotted,arcangle=20]{71}{73} \ncline{71}{81} \ncline{72}{82}  \ncline{73}{83}
\ncline{81}{82} \ncline{82}{83} \ncarc[linestyle=dotted,arcangle=20]{81}{83} \ncline{81}{91} \ncline{82}{92}  \ncline{83}{93}
\ncline{91}{92} \ncline{92}{93} \ncarc[linestyle=dotted,arcangle=20]{91}{93} \ncline{91}{101} \ncline{92}{102}  \ncline{93}{103}
\ncline{101}{102}  \ncline{102}{103} \ncarc[linestyle=dotted,arcangle=20]{101}{103}

\ncarc[linestyle=dashed]{11}{101}  \ncarc[linestyle=dashed]{12}{102}  \ncarc[linestyle=dashed]{13}{103}

\end{pspicture}
\caption{\label{f2}\footnotesize The minimum paired dominating set
(bold vertices) of $G_{10,3}$}
\end{center}
\end{figure}
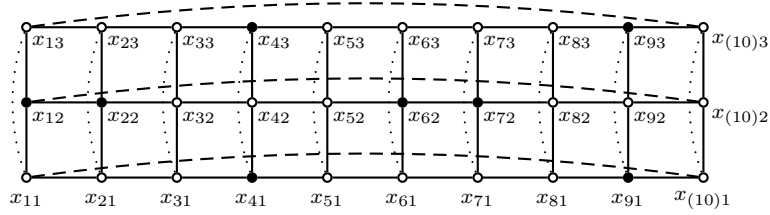

To obtain the upper bounds of $\gamma_{t}(G_{n,3})$ and
$\gamma_{p}(G_{n,3})$, we set
 $$
 D=\{x_{i2}:
 i\equiv\,1,2 \,({\rm mod}\,5)\} \cup \{x_{j1}, x_{j3}: j\equiv\,4
 \,({\rm mod}\,5)\},
 $$
where $1\le i\le n$. See Figure~\ref{f2}, where $D$ consists of bold
vertices.

If $n\not\equiv \,3 \,({\rm mod}\,5)$, then $D$ is a total
dominating set and $\gamma_{t}(G_{n,3})\leq |D|=\lceil
\frac{4n}{5}\rceil$.

If $n\equiv \,3 \,({\rm mod}\,5)$, then $D\cup \{x_{n2}\}$ is a
total dominating set and $\gamma_{t}(G_{n,3})\leq |D|+1=\lceil
\frac{4n}{5}\rceil$.

Combining these facts with (\ref{e3.4}), we have that
$\gamma_{t}(G_{n,3})=\lceil \frac{4n}{5}\rceil$.

If $n\equiv \,0,2,4 \,({\rm mod}\,5)$, then $D$ is a paired
dominating set and $\gamma_{p}(G_{n,3})\leq |D|=\lceil
\frac{4n}{5}\rceil$.

If $n\equiv \,1\,({\rm mod}\,5)$, then $D\cup \{x_{n1}\}$ is a
paired dominating set and $\gamma_{p}(G_{n,3})\leq |D|+1=\lceil
\frac{4n}{5}\rceil+1$.

If $n\equiv \,3 \,({\rm mod}\,5)$, then $D\cup \{x_{n1},x_{n2}\}$ is
a paired dominating set and $\gamma_{p}(G_{n,3})\leq |D|+2=\lceil
\frac{4n}{5}\rceil+1$.

Since $\gamma_p(G_{n,3})\ge \gamma_t(G_{n,3})$ and
$\gamma_p(G_{n,3})$ is even,
$\gamma_p(G_{n,3})=\lceil\frac{4n}{5}\rceil$ if $n\equiv \,0,2,4
\,({\rm mod}\,5)$, and
$\gamma_p(G_{n,3})=\lceil\frac{4n}{5}\rceil+1$ if $n\equiv\,1,3
\,({\rm mod}\,5)$.

The theorem follows.
\end{pf}


\section{Total and paired domination number of $G_{n,4}$}

In this section, we determine the exact values of
$\gamma_{t}(G_{n,4})$ and $\gamma_{p}(G_{n,4})$, the latter has been
announced by Bre\v{s}ar, Henning and Rall~\cite{bhr05}, but without
proofs.

\begin{lem}\label{lem4.2}
$\gamma_{p}(G_{n,4})=\gamma_{t}(G_{n,4})=n+1$ for $n\equiv
1,3~\,({\rm mod}~4)$.
\end{lem}

\begin{pf}
For $n\equiv 1~\,({\rm mod}~4)$, let
 $$
 D=\{x_{i1},x_{i2},x_{(i+2)3},x_{(i+2)4}:\ i\equiv 1\,({\rm
mod}\, 4), i\ne n\} \cup \{x_{n1},x_{n2}\}.
$$
Then $D$ is a paired dominating set of $G_{n,4}$ with cardinality
$n+1$. For $n\equiv 3~\,({\rm mod}~4)$,
$D=\{x_{i1},x_{i2},x_{(i+2)3},x_{(i+2)4}:~i\equiv 1~\,({\rm
mod}~4)\}$
is a paired dominating set of $G_{n,4}$ with cardinality $n+1$.
Thus, $\gamma_{t}(G_{n,4})\le \gamma_{p}(G_{n,4})\le n+1$ for
$n\equiv 1,3~\,({\rm mod}~4)$.

By Lemma~\ref{lem2.2}, $\gamma_{t}(G_{n,4})\ge \frac{4n}{4}=n$. Now,
we prove $\gamma_{t}(G_{n,4})\ge n+1$. Suppose to the contrary that
$\gamma_{t}(G_{n,4})=n$. By Lemma~\ref{lem2.2}, $G_{n,4}$ has an
efficient total dominating set $D'$. By Lemma~\ref{lem2.1}, $|D'|=n$
is even, a contradiction. Therefore $\gamma_{t}(G_{n,4})>n$, and
hence $\gamma_{p}(G_{n,4})=\gamma_{t}(G_{n,4})=n+1$.
\end{pf}

\begin{lem}\label{lem4.3}
$\gamma_{t}(G_{n,4})\le \gamma_{p}(G_{n,4})\le n+2$ for $n\equiv
2~\,({\rm mod}~4)$.
\end{lem}

\begin{pf}
Let
 $$
 D=\{x_{i1},x_{i2},x_{(i+2)3},x_{(i+2)4}:\ i\equiv 1\,({\rm
 mod}\,4), i\le n-2\} \cup \{x_{(n-1)1},x_{(n-1)2},x_{n1},x_{n2}\}.
 $$
Then $D$ is a paired dominating set of $G_{n,4}$ with cardinality
$n+2$. Thus, $\gamma_{t}(G_{n,4})\le \gamma_{p}(G_{n,4})\le n+2$.
\end{pf}

\vskip6pt

To prove $\gamma_{t}(G_{n,4})\ge n+2$ for $n\equiv 2~\,({\rm
mod}~4)$, we need the following notations and two lemmas. Let
$H_i^j=Y_{i}\cup Y_{i+1} \cup\ldots \cup Y_{i+j-1}$, and let $G_i^j$
be the graph obtained from $G_{n,4}-H_i^j$ by adding the edge-set
$\{x_{(i-1)k}x_{(i+j)k}:\ 1\le k\le 4\}$, where the subscripts are
modulo $n$. Clearly, $G_i^j\cong G_{n-j,4}$.
%

\begin{lem}\label{lem4.4}
Let $D$ be a total dominating set of $G_{n,4}$. Then $|D\cap
H_i^4|\ge 4$ for any $i$ with $1\le i\le n$. Moreover, if there
exists some $i$ with $1\le i\le n$ such that $|N(v)\cap D|=1$ for
any vertex $v$ in $H_i^4$, then $D'=D\setminus (D\cap H_i^4)$ is a
total dominating set of $G_i^4$.
\end{lem}

\begin{pf}
Without loss of generality, assume $i=2$. It can be easy verified to
dominate 8 vertices in $Y_3\cup Y_4$, at least $4$ vertices are
needed, and hence $|D\cap H_2^4|\ge 4$.

We now show the second assertion. 
Suppose to the contrary that $D'$ is not a total dominating set of
$G_2^4$. Then there is a vertex $u$ in $Y_1\cup Y_6$ such that it is
not dominated by $D'$, that is, $N_{G_2^4}(u)\cap D'=\emptyset$.
Without loss of generality assume $u=x_{11}$. Then $x_{21}\in D$ and
$x_{61}\notin D$. Also $x_{41}\notin D$ since $|N(x_{31})\cap D|=1$.

Since $x_{33}$ should be dominated by $D$ and $|N(x_{33})\cap D|=1$,
only one of $x_{32}$, $x_{34}$, $x_{23}$, and $x_{43}$ belongs to
$D$. If $x_{32}\in D$ or $x_{34}\in D$, then $|N(x_{31})\cap D|\ge
2$, a contradiction. If $x_{23}\in D$, then $|N(x_{22})\cap D|\ge
2$, a contradiction. Thus, $x_{43}\in D$. Since $x_{51}$ should be
dominated by $D$, $x_{52}\in D$ or $x_{54}\in D$. But then
$|N(x_{53})\cap D|\ge 2$, a contradiction. Thus, $D'=D\setminus
(D\cap H_2^4)$ is a total dominating set of $G_i^4$.
\end{pf}

\begin{lem}\label{lem4.5}
Let $D$ be a total dominating set of $G_{n,4}$. If $x_{ij}$ is
dominated by two vertices $u,v\in D$, then there exists a vertex $w$
in $H_{i-1}^2$ or $H_i^2$ such that $|N(w)\cap D|\ge 2$.
\end{lem}

\begin{pf}
Without loss of generality, let $i=j=2$. If $u,v\in Y_2$, then
assume $u=x_{21}$, $v=x_{23}$ and, hence, $|N(x_{24})\cap D|\ge 2$.

If one of $u$ and $v$ is in $Y_2$ and another is in $Y_1\cup Y_3$, then
without loss of generality assume $u=x_{21}\in Y_2$ and $v=x_{32}\in Y_3$.
And then $|N(x_{31})\cap D|\ge 2$.

If one of $u$ and $v$ is in $Y_1$ and another is in $Y_3$, then
without loss of generality assume $u=x_{12}\in Y_2$ and $v=x_{32}\in
Y_3$. Since $x_{24}$ should be dominated by $D$, let $s\in
N(x_{24})\cap D$. It is clearly that $N(s)\cap N(u)\ne \emptyset$ or
$N(s)\cap N(v)\ne \emptyset$, which implies that there exists a
vertex $w\notin \{u,v\}$ in $H_{1}^2\cup H_{2}^2$ such that
$|N(w)\cap D|\ge 2$.
\end{pf}

\begin{lem}\label{lem4.6}
$\gamma_{t}(G_{n,4})=\gamma_{p}(G_{n,4})= n+2$ for $n\equiv
2~\,({\rm mod}~4)$.
\end{lem}

\begin{pf}
By Lemma~\ref{lem4.3}, we only need to show $\gamma_{t}(G_{n,4})\ge
n+2$. To this end, let $n=4k+2$. We proceed by induction on $k\ge
1$. It is easy to verify that $\gamma_{t}(G_{6,4})=8$ and
$\gamma_{t}(G_{10,4})=12$. The conclusion is true for $k=1,2$.
Assume that the induction hypothesis is true for $k-1$ with $k\ge
3$.

Let $D$ be a minimum total dominating set of $G_{n,4}$, where
$n=4k+2$ for $k\ge 3$. Assume to the contrary that $|D|\le n+1$.
Since any vertex $u$ can dominate at most 4 vertices in $G_{n,4}$
and $|V(G_{n,4})|=4n$, there are at most four vertices such that
each of them is dominated by at least two vertices in $D$.

We now prove that there exists some $i\in\{1,2,\ldots, n\}$ such
that $|N(v)\cap D|=1$ for any vertex $v\in H_i^4$. There is nothing
to do if there are at most three vertices such that each of them is
dominated by at least two vertices since $n\ge 14$. Now, assume
there are exactly four vertices such that each of them is dominated
by at least two vertices. By Lemma~\ref{lem4.5}, there exists two
integers $s$ and $t$ with $1\le s,t\le n$ such that two of the four
vertices are in $H_s^2$ and the other two are in $H_t^2$. Therefore,
there exists an integer $i$ with $1\le i\le n$ such that for any
vertex $v\in Y_i$, $|N(v)\cap D|=1$ since $n\ge 14$.

By Lemma~\ref{lem4.4}, $|D\cap H_i^4|\ge 4$ and $D'=D\setminus
(D\cap H_i^4)$ is a total dominating set of $G_i^4\cong G_{n-4,4}$.
By the inductive hypothesis, $|D'|\ge \gamma_t(G_{n-4,4})\ge n-2$.
It follows that
 $$
 n+1\ge |D|=|D\cap H_i^4|+|D'|\ge 4+n-2=n+2,
 $$
a contradiction, which implies that $\gamma_{t}(G_{n,4})=|D|\ge
n+2$. By the induction principle, the lemma follows.
\end{pf}

\vskip6pt

We state the above results as the following theorem.

\begin{thm}\label{thm4.1} For any integer $n\ge 3$,
 $$
\gamma_{t}(G_{n,4})=\gamma_{p}(G_{n,4})=\left\{
 \begin{array}{ll}
n, & {\rm if}\ n\equiv 0\,({\rm mod}\,4);\\
n+1,& {\rm if}\ n\equiv 1,3\,({\rm mod}\,4);\\
n+2,& {\rm if}\ n\equiv 2\,({\rm mod}\,4).
 \end{array}\right.$$
\end{thm}

\section{Upper bounds of $\gamma_{p}(G_{n,m})$ for $n, m\ge 5$}

The values of $\gamma_{t}(G_{n,m})$ and $\gamma_{p}(G_{n,m})$ for
$m\in\{3,4\}$ have been determined in the above sections, but their
values for $m\ge 5$ have been not determined yet. In this section,
we present their upper bounds. Since $\gamma_t(G)\le \gamma_p(G)$
for any graph $G$ without isolated vertices, we establish upper
bounds only for $\gamma_{p}(G_{n,m})$ if we can not obtain a smaller
upper bound of $\gamma_t(G_{n,m})$ than that of
$\gamma_{p}(G_{n,m})$.

\begin{lem}\label{lem5.1}
$\gamma_t(G_{n,m})\le \gamma_t(G_{n+1,m})$ and $\gamma_p(G_{n,m})\le
\gamma_p(G_{n+1,m})$.
\end{lem}

\begin{pf}
Let $D$ be a minimum paired (total) dominating set of $G_{n+1,m}$.

If $D\cap Y_{n+1}=\emptyset$, then $D$ is also a paired (total)
dominating set of $G_{n,m}$, and hence $\gamma_{p}(G_{n,m})\leq |D|$
($\gamma_{t}(G_{n,m})\leq |D|$).

Assume $D\cap Y_{n+1}\ne\emptyset$ below. Let $A=\{j|\ x_{(n+1)j}\in
D\}$ and $B=\{j|\ x_{nj}\in D\}$. Then $D'=(D\setminus Y_{n+1})\cup
\{x_{(n-1)j}|\ j\in A\cap B\} \cup \{x_{nj}|\ j\in A\setminus B\}$
is a total dominating set of $G_{n,m}$ and $|D'|\leq |D|$. Therefore
$\gamma_{t}(G_{n,m})\leq \gamma_t(G_{n+1,m})$.

The vertex set $D'$  may not be a paired dominating set of
$G_{n,m}$, that means, the induced subgraph $G$ by $D'$ in $G_{n,m}$
may contains odd connected components. Let $p$ be the number of odd
connected components in $G$. It is clear that $|D'|\le |D|-p$ by the
construction of $D'$ from $D$. Therefore, we can obtain $D''$ by
adding at most $p$ vertices to $D'$ such that the induced subgraph
by $D''$ in $G_{n,m}$ does not contain odd connected components.
Then $D''$ is a paired dominating set of $G_{n,m}$, and hence
$\gamma_{p}(G_{n,m})\leq |D''|\leq |D|$.
\end{pf}

\begin{thm}\label{thm5.1}
$\gamma_p(G_{n,m})\le
4\lceil\frac{n}{4}\rceil\lceil\frac{m}{4}\rceil$.
\end{thm}

\begin{pf}
Let $n=4a-i$ and $m=4b-j$ where $0\le i,j\le 3$. By
Lemma~\ref{lem2.3}, $\gamma_p(G_{4a,4b})=4ab
=4\lceil\frac{n}{4}\rceil\lceil\frac{m}{4}\rceil$. By
Lemma~\ref{lem5.1}, $\gamma_p(G_{m,n})\le \gamma_p(G_{4a,4b})
=4\lceil\frac{n}{4}\rceil\lceil\frac{m}{4}\rceil$.
\end{pf}

\vskip6pt

For $n,m\ge 5$, let $m\equiv \,a~({\rm mod}\,4)$ and $n\equiv
\,b~({\rm mod}\,4)$ where $0 \le a,b\le 3$. We will establish some
better bounds of $\gamma_{t}(G_{n,m})$ and $\gamma_{p}(G_{n,m})$
than those in Theorem~\ref{thm5.1} for some special $a$ and $b$.
Without loss of generality, we can assume $b\ge a$ since
$G_{n,m}\cong G_{m,n}$. Let
 $$
 D_e=\{x_{ij},x_{i(j+1)},x_{(i+2)(j+2)},x_{(i+2)(j+3)}: ~i,j\equiv
\,1~({\rm mod}\,4)\},
$$ where $1\le i\le n-2$, $1\le j\le m-2$, and $n,m\ge 5$.

\begin{thm}
$\gamma_{p}(G_{n,m})\le \frac{(n+1)m}{4}$ for $m\equiv \,0~({\rm
mod}\,4)$ and $n\equiv \,1~({\rm mod}\,4)$.
\end{thm}

\begin{pf}
Let $D=D_e\cup \{x_{nj},x_{n(j+1)}: ~j\equiv \,1~({\rm mod}\,4)\}$,
where $1\le j\le m-2$. Then, it is easy to see that $D$ is a paired
dominating set of $G_{n,m}$ with cardinality $\frac{(n+1)m}{4}$.
Thus, $\gamma_{p}(G_{n,m})\le \frac{(n+1)m}{4}$.
\end{pf}

\begin{thm}\label{thm5.3}
$\gamma_{t}(G_{n,m})\le  \frac{(n+1)(m+1)}{4}$ and
$\gamma_{p}(G_{n,m})\le \frac{(n+1)(m+1)}{4}+1$ for $m,n\equiv
\,1~({\rm mod}\,4)$.
\end{thm}

\begin{pf}
Let $D=D_e\cup\{x_{nj},x_{n(j+1)},x_{(i+1)(m-1)},x_{(i+2)m}:
~i,j\equiv \,1~({\rm mod}\,4)\}\cup \{x_{nm}\}$, where $1\le i\le
n-2$ and $1\le j\le m-2$. Then, it is easy to see that $D$ is a
total dominating set of $G_{n,m}$ with cardinality
$\frac{(n+1)(m+1)}{4}$, and $D\cup \{x_{n(m-1)}\}$ is a paired
dominating set of $G_{n,m}$ with cardinality
$\frac{(n+1)(m+1)}{4}+1$. Thus, $\gamma_{t}(G_{n,m})\le
\frac{(n+1)(m+1)}{4}$ and $\gamma_{p}(G_{n,m})\le
\frac{(n+1)(m+1)}{4}+1$.
\end{pf}

\begin{thm}\label{thm5.4}
$\gamma_{t}(G_{n,m})\le \frac{(n+1)(m+1)}{4}-3$ and
$\gamma_{p}(G_{n,m})\le \frac{(n+1)(m+1)}{4}-2$ for $m\equiv
\,1~({\rm mod}\,4)$ and $n\equiv \,3~({\rm mod}\,4)$.
\end{thm}

\begin{pf}
Let $D=(D_e\cup\{x_{(i+1)(m-1)},x_{(i+2)m} :~i\equiv \,1~({\rm
mod}\,4)\})\setminus \{x_{n(m-2)},x_{nm}\}$, where $1\le i\le n-2$.
Then, $D$ is a paired dominating set of $G_{n,m}$ with cardinality
$\frac{(n+1)(m+1)}{4}-2$, and $D\setminus\{x_{2(m-1)}\}$ is a total
dominating set of $G_{n,m}$ with cardinality
$\frac{(n+1)(m+1)}{4}-3$. Thus, $\gamma_{t}(G_{n,m})\le
\frac{(n+1)(m+1)}{4}-3$ and $\gamma_{p}(G_{n,m})\le
\frac{(n+1)(m+1)}{4}-2$.
\end{pf}

\begin{cor}
$\gamma_{t}(G_{n,m})\le \frac{(n+2)(m+1)}{4}-3$ and
$\gamma_{p}(G_{n,m})\le \frac{(n+2)(m+1)}{4}-2$ for $m\equiv
\,1~({\rm mod}\,4)$ and $n\equiv \,2~({\rm mod}\,4)$.
\end{cor}

\begin{pf}
By Lemma~\ref{lem5.1}, $\gamma_t(G_{n,m})\le \gamma_t(G_{n+1,m})$
and $\gamma_p(G_{n,m})\le \gamma_p(G_{n+1,m})$. The corollary
follows from Theorem~\ref{thm5.4}.
\end{pf}

\begin{thm}\label{thm5.5}
$\gamma_{p}(G_{n,m})\le \frac{(n+2)(m+2)}{4}-6$ for $m,n\equiv
\,2~({\rm mod}\,4)$.
\end{thm}
\begin{pf}
Let $D=(D_e\cup\{x_{i(m-2)},x_{i(m-1)},x_{(i+2)(m-1)},x_{(i+2)m}:\
i\equiv \,1~({\rm mod}\,4)\}\cup \{x_{(n-1)j},$
$x_{(n-1)(j+1)},x_{n(j+2)},x_{n(j+3)}:\ j\equiv \,1~({\rm
mod}\,4)\}\cup \{x_{n(m-1)}\}) \setminus \{x_{1(m-2)},$
$x_{1(m-1)},x_{n(m-3)}\}$, where $1\le i\le n-2$ and $1\le j\le
m-2$. Then $D$ is a paired dominating set of $G_{n,m}$ with
cardinality $\frac{(n+2)(m+2)}{4}-6$. Thus, $\gamma_{p}(G_{n,m})\le
\frac{(n+2)(m+2)}{4}-6$.
\end{pf}

\end{document}